\documentclass[12pt]{amsart}
\usepackage{amsmath,mathrsfs}

\usepackage[ps2pdf,colorlinks=true,urlcolor=blue,
citecolor=red,linkcolor=blue,linktocpage,pdfpagelabels,bookmarksnumbered,bookmarksopen]{hyperref}
\usepackage[english]{babel}

\usepackage[left=2.8cm,right=2.8cm,top=2.5cm,bottom=2.5cm]{geometry}

\numberwithin{equation}{section}

\newcommand{\R}{{\mathbb R}}
\newcommand{\T}{{\mathbb T}}

\newcommand{\eps}{\varepsilon}

\newcommand{\ws}[1]{|d#1|}

\renewcommand{\theta}{\vartheta}

\numberwithin{equation}{section}
\newtheorem{theorem}{Theorem}[section]
\newtheorem{proposition}[theorem]{Proposition}

\newtheorem{corollary}[theorem]{Corollary}

\theoremstyle{definition}

\title{On Ekeland's variational principle}

\author{Marco Squassina}
\address{Dipartimento di Informatica
\newline\indent
Universit\`a degli Studi di Verona
\newline\indent
C\'a Vignal 2, Strada Le Grazie 15, I-37134 Verona, Italy}
\email{marco.squassina@univr.it}

\thanks{Research supported by PRIN: {\em Metodi Variazionali e Topologici
nello Studio di Fenomeni non Lineari}}

\begin{document}
	



\begin{abstract}
For proper lower semi-continuous functionals bounded below which do not increase upon polarization, an improved version of Ekeland's 
variational principle can be formulated in Banach spaces, which provides almost symmetric points. 
\end{abstract}

\maketitle

\section{Introduction}
In the study of nonconvex minimization problems, the idea of looking for a special minimizing sequence with good properties
in order to guarantee the convergence towards a minimizer goes back to the work of Hilbert and Lebesgue \cite{hilbert,lebesgue}.
In this direction one of the main contributions of the last decades in the calculus of variations was surely
provided by Ekeland's variational principle for lower semi-continuous 
functionals on metric spaces, discovered in 1972 \cite{ekeland0,ekeland1}.
Since then, it has found a multitude of applications in different fields of nonlinear analysis and turned out to be fruitful in simplifying 
and unifying the proofs of already known results. We refer the reader to the survey \cite{ekeland2} and 
to the monograph \cite{aubin-ekeland} for a discussion on a 
broad range of applications, including optimization, control and geometry of Banach spaces.
The aim of the present note is that of showing that within the abstract symmetrization framework proposed by Van Schaftingen
in the nice paper \cite{vansch}, 
under the assumption that the functional does not increase by polarization, the conclusion of the principle can be
enriched with the useful information that the existing almost critical point is almost symmetric as well. Furthermore,
in the context of the critical point theory for nonsmooth functionals originally developed in \cite{dm},
the result allows to detect a Palais-Smale sequence $(u_h)$ in the sense of weak slope, which becomes more and more
symmetric, as $h\to\infty$, yielding a symmetric minimum point, provided that some compactness is available. 
The additional feature often gives rise to a compactifying effect through suitable compact embeddings of spaces of symmetric functions.
For minimax critical values of $C^1$ functionals, a similar result has been obtained in \cite{vansch} under the assumption that
the functional enjoys a mountain pass geometry. Recently the whole machinery has been extended by the author in \cite{squass} to a class of
lower semi-continuous functionals in the framework of~\cite{dm}. 
Deriving a symmetric version of Ekeland's variational principle in Banach spaces
is easier than obtaining the results of \cite{vansch,squass}, since handling a minimization sequence is of course simpler 
than managing a minimaxing sequence. On the other hand, due to the great impact of Ekeland's principle
in the mathematical literature over the last three decades, the author believes that 
highlighting the precise statements could reveal useful for various applications.
Let us now come to the formulation of the results.
Let $X$ and $V$ be two Banach spaces and $S\subseteq X$. We shall consider two maps $*:S\to S$, $u\mapsto u^*$, 
the symmetrization map, and $h:S\times {\mathcal H}_*\to S$, $(u,H)\mapsto u^H$, the polarization map, ${\mathcal H}_*$ 
being a path-connected topological space. We assume, according to \cite[Section 2.4]{vansch}, that the following hold:
\begin{enumerate}
 \item $X$ is continuously embedded in $V$;
 \item $h$ is a continuous mapping;
\item for each $u\in S$ and $H\in {\mathcal H}_*$ it holds $(u^*)^H=(u^H)^*=u^*$ and $u^{HH}=u^H$;
\item there exists $(H_m)\subset {\mathcal H}_*$ such that, for $u\in S$, $u^{H_1\cdots H_m}$ converges
to $u^*$ in $V$;
\item for every $u,v\in S$ and $H\in {\mathcal H}_*$ it holds
$\|u^H-v^H\|_V\leq \|u-v\|_V$.
\end{enumerate}
Moreover, the mapping $*:S\to V$ can be extended to $*:X\to V$ by 
setting $u^*:=(\Theta(u))^*$ for every $u\in X$, where $\Theta:(X,\|\cdot\|_V)\to (S,\|\cdot\|_V)$ is
a Lipschitz function, of Lipschitz constant $C_\Theta>0$, such that $\Theta|_{S}={\rm Id}|_{S}$. We refer to 
\cite[Section 2.4]{vansch} for some examples of concrete situations suitable in applications 
to partial differential equations. 
\vskip5pt
\noindent
We recall~\cite[Corollary 3.1]{vansch} a useful result on approximation of symmetrizations.

\begin{proposition}
	\label{mapJvS}
For all $\rho>0$ there exists a continuous mapping $\T_\rho:S\to S$ such that $\T_\rho u$ 
is built via iterated polarizations and $\|\T_\rho u-u^*\|_V<\rho$, for all $u\in S$.
\end{proposition}

\noindent
In the above framework, here is the main result.

\begin{theorem}
	\label{mainthm}
Assume that $f:X\to\R\cup\{+\infty\}$ is a proper and lower 
semi-continuous functional bounded from below such that
\begin{equation}
	\label{assumptionpol}
\text{$f(u^H)\leq f(u)$\,\,\quad for all $u\in S$ and $H\in {\mathcal H}_*$}.
\end{equation}
Let $u\in S$, $\rho>0$ and $\sigma>0$ with
$$
f(u)<\inf_X f+\rho\sigma.
$$
Then there exists $v\in X$ such that
\begin{enumerate}
\item[(a)] $\|v-v^*\|_V\leq C\rho$; 
\item[(b)] $\|v-u\|\leq \rho+\|\T_\rho u-u\|;$ 
\item[(c)]   $f(v)\leq f(u)$;   
\item[(d)]   $f(w)\geq f(v)-\sigma \|w-v\|,$\quad\text{for all $w\in X,$}
\end{enumerate}
for some positive constant $C$ depending only upon $V,X$ and $\Theta$.
\end{theorem}

\noindent
Denoting by $|df|(u)$ the weak slope \cite{dm} of $f$ at $u$ (it is $|df|(u)=\|df(u)\|_{X'}$ if $f$ is $C^1$),
we say that $(u_j)\subset X$ is a {\em symmetric Palais-Smale sequence} 
at level $c\in\R$ ($(SPS)_c$-sequence) if $\ws{f}(u_j) \to 0$, $f(u_j) \to c$ and, in addition, 
$\|u_j-u_j^*\|_V\to 0$ as $j\to\infty$. 
We say that $f$ satisfies the {\em symmetric Palais-Smale condition} at level $c$ ($(SPS)_c$
in short), if any $(SPS)_c$ sequence has a subsequence converging in $X$.
In this context, we have the following

\begin{corollary}
		\label{maincor}
Assume that $f:X\to\R\cup\{+\infty\}$ is a proper and lower 
semi-continuous functional bounded from below which 
satisfies~\eqref{assumptionpol}. Moreover, assume that
for all $u\in X$ there exists $\xi\in S$ with $f(\xi)\leq f(u)$.
Then, for any $\eps>0$, there exists $v\in X$ such that
\begin{enumerate}
\item[(a)] $\|v-v^*\|_V\leq C\eps$; 
\item[(b)]   $f(v)<{\mathcal M}+\eps^2$,\quad ${\mathcal M}:=\inf f$;   
\item[(c)]   $|df|(v)\leq \eps$,
\end{enumerate}
for some $C>0$.
In particular, $f$ has a $(SPS)_{\mathcal M}$-sequence. If $f$ satisfies $(SPS)_{\mathcal M}$ 
then $f$ admits a critical point (for the weak slope) $z\in X$  with $f(z)={\mathcal M}$ and $z=z^*$. 
\end{corollary}

This statement sounds particularly useful for applications to PDEs, via 
the additional control (a) that often gives rise to compactifying effects 
(see, for instance, the arguments of \cite[pages 479 and 480]{vansch}) via suitable compact embeddings 
of spaces of symmetric functions (we refer, for instance, to \cite[Section I.1.5]{willem}). We note that
the assumption that, for all $u\in X$, there exists an element $\xi\in S$ such that $f(\xi)\leq f(u)$
is satisfied in many typical concrete situations, like when $X$ is a Sobolev space, $S$ is the cone of its positive
functions and the functional satisfies $f(|u|)\leq f(u)$, for all $u\in X$.  
\vskip4pt
\noindent
Finally, let ${\mathcal B}^*_{{\mathcal H}_*}$ the set of $\varphi\in X^*$ such that $\|\varphi\|\leq 1$,
$\langle\varphi,u\rangle\leq \langle\varphi,u^H\rangle$
for $u\in S$ and $H\in {\mathcal H}_*$ and for any $u\in X$ there is $\xi\in S$ 
with $\langle\varphi,u\rangle\leq \langle\varphi,\xi\rangle$, being $X^*$ the dual of $X$. 
In the spirit of \cite[Corollary 2.4]{ekeland1}, Corollary~\ref{maincor} 
yields the following density result.
\begin{corollary}
		\label{corgeom}
Assume that $f:X\to\R$ is a G\^ateaux differentiable function satisfying the assumptions
of Corollary~\ref{maincor}. Moreover, let $\alpha>0$ and $\beta\in\R$ be such that
$$
f(v)\geq \alpha \|v\|+\beta,\quad\text{for all $v\in X$}.
$$
Then, if ${\mathcal S}:=\{v\in X:\|v-v^*\|_V\leq 1\}$, the set $df({\mathcal S})\subset X^*$ 
is dense in $\alpha {\mathcal B}^*_{{\mathcal H}_*}$.
\end{corollary}

\section{Proofs}

\subsection{Proof of Theorem~\ref{mainthm}}
Let $u\in S$, $\rho>0$ and $\sigma>0$ be such that $f(u)<\inf f+\rho\sigma$.
If $\T_\rho:S\to S$ is the continuous mapping of Proposition~\ref{mapJvS},
we set $\tilde u:=\T_\rho u\in S$. Then, by construction we have $\|\tilde u-u^*\|_V<\rho$
and, in light of~\eqref{assumptionpol} and the property that $\tilde u$ is built 
from $u$ through iterated polarizations, we obtain
$$
f(\tilde u)<\inf_X f+\rho\sigma.
$$
By Ekeland's variational principle (cf.~\cite{ekeland0,ekeland1,ekeland2}), there exists
an element $v\in X$ such that
$$
f(v)\leq f(\tilde u),\quad\,\,
\|v-\tilde u\|\leq \rho,\quad\,\,
f(w)\geq f(v)-\sigma \|w-v\|,\,\,\,\text{for all $w\in X$}.
$$
Hence (d) holds and, since $f(v)\leq f(\tilde u)\leq f(u)$, conclusion (c) follows as well.
In the abstract symmetrization framework, it is readily seen that
$\|u^*-v^*\|_V\leq C_\Theta\|u-v\|_V$ for all $u,v\in X$. Then, if $K>0$ is the 
continuity constant of the injection $X\hookrightarrow V$, it follows
\begin{align*}
\|v-v^*\|_V &\leq \|v-\tilde u\|_V+\|\tilde u-u^*\|_V+\|u^*-v^*\|_V \\
&\leq K(C_\Theta+1)\|v-\tilde u\|+\|\tilde u-u^*\|_V  \leq (K(C_\Theta+1)+1)\rho,
\end{align*}
where we used the fact that $u^*=\tilde u^*$, in light of (3) of the
abstract framework and, again, by the way $\tilde u$ is built from $u$. 
Then, also conclusion (a) holds true. Finally, we have
\begin{equation}
\|v-u\|\leq \|v-\tilde u\|+ \|\tilde u-u\|\leq \rho+\|\T_\rho u-u\|,
\end{equation}
yielding (b). This concludes the proof of the theorem. \qed
\vskip2pt

\subsection{Proof of Corollary~\ref{maincor}}

Given $\eps>0$, let $u\in X$ be such that $f(u)<{\mathcal M}+\eps^2$. By assumption, 
we can find an element $\hat u\in S$ such that 
\begin{equation}
	\label{firstineqq}
f(\hat u)<{\mathcal M} +\eps^2.
\end{equation}
Then, we are allowed to apply Theorem~\ref{mainthm} with $\sigma=\rho=\eps$ and get an element $v\in X$
such that $\|v-v^*\|_V\leq C\eps$ for some positive constant $C$ depending only upon $V,X$ and $\Theta$, $f(v)\leq f(\hat u)$ and 
$f(w)\geq f(v)-\eps \|w-v\|$, for all $w\in X$. Whence, by taking into account inequality \eqref{firstineqq}, 
conclusions (a) and (b) of the corollary hold true. Moreover, since 
$$
\limsup_{w\to v}\frac{f(v)-f(w)}{\|v-w\|}\leq \eps,
$$
by the definition of strong slope $|\nabla f|(u)$ \cite{dm} it follows $|\nabla f|(u)\leq\eps$.
Hence, since the weak slope satisfies $|df|(u)\leq |\nabla f|(u)$ \cite{dm}, assertion (c) immediately follows. 
Choosing now a sequence $(\eps_j)\subset\R^+$ with $\eps_j\to 0$ as $j\to\infty$, 
by definition one finds a $(SPS)_{\mathcal M}$-sequence $(v_j)\subset X$. If $f$ satisfies $(SPS)_{\mathcal M}$,
there exists a subsequence, that we still denote by $(v_j)$, which converges to some $z$ in $X$. Hence, 
via lower semi-continuity, we get
$$
{\mathcal M}\leq f(z)\leq \liminf_{j\to\infty} f(v_j)={\mathcal M}.
$$
Since $|df|(v_j)\to 0$ and $f(v_j)\to f(z)={\mathcal M}$ as $j\to\infty$, by means of \cite[Proposition 2.6]{dm}, it follows
that $|df|(z)\leq\liminf_j |df|(v_j)=0$. Notice that, since $\|v_j-v_j^*\|_V\to 0$ as $j\to\infty$, letting $j\to\infty$ 
into the inequality
\begin{equation*}
\|z-z^*\|_V \leq \|z-v_j\|_V+\|v_j-v_j^*\|_V+\|v_j^*-z^*\|_V 
 \leq (C_\Theta+1)K\|v_j-z\|+\|v_j-v_j^*\|_V,  
\end{equation*}
yields $z=z^*$, as desired. This concludes the proof of the corollary. \qed

\bigskip

\end{document}